\documentclass[12pt]{amsart}
\usepackage{amsmath, amsthm, amsfonts, amssymb}

\newtheorem{theorem}{Theorem}[section]
\newtheorem{lemma}[theorem]{Lemma}

\theoremstyle{definition}

\newtheorem{example}[theorem]{Example}

\theoremstyle{remark}
\newtheorem{remark}[theorem]{Remark}
\numberwithin{equation}{section}

\newcommand{\mC}{\ensuremath{\mathbb{C}}}
\newcommand{\mD}{\ensuremath{\mathbb{D}}}

\newcommand{\mN}{\ensuremath{\mathbb{N}}}

 \allowdisplaybreaks
 
\begin{document}

\title{Sharing of a Set of Meromorphic Functions and Montel's Theorem}

\author[K. S. Charak]{ Kuldeep Singh Charak}
\address{Department of Mathematics, University of Jammu,
Jammu-180 006, India}
\email{kscharak7@rediffmail.com}

\author[V. Singh]{Virender Singh}
\address{Department of Mathematics, University of Jammu,
Jammu-180 006, India}
\email{virendersingh2323@gmail.com }

\begin{abstract}
In this paper we prove the result:
Let $\mathcal{F}$ be a family of meromorphic functions on a domain $\Omega$ such that every pair of members of
$\mathcal{F}$  shares a set $S:=\left\{\psi_1(z), \psi_2(z), \psi_3(z) \right\}$ in $\Omega$, where $\psi_j(z), \ j=1,2,3$ is meromorphic in $\Omega.$ 
If for every $f\in \mathcal{F}$, $f(z_0)\neq \psi_i (z_0)$ whenever $\psi_i(z_0)=\psi_j(z_0)$ for $i,j\in \left\{1,2,3 \right\}(i\neq j)$ and $z_0\in \Omega ,$ 
then $\mathcal{F}$ is normal in $\Omega$. This result generalizes a result of M.Fang and W.Hong [{\it Some results on normal family of meromorphic functions, Bull. Malays. Math. Sci. Soc. (2)23 (2000),143-151,}] and in particular, it generalizes the most celebrated theorem of Montel-the Montel's theorem.

\end{abstract}

\renewcommand{\thefootnote}{\fnsymbol{footnote}}
\footnotetext{The research work of the second author is supported by the CSIR India.}

\maketitle

\section{\textbf{Introduction and Main Results}}
A family $\mathcal F$ of meromorphic functions defined on a domain $\Omega \subseteq \overline {\mC}$ is said to be normal in $\Omega$ if every sequence of elements of $\mathcal F$  contains a subsequence which converges locally uniformly in $\Omega$ with respect to the spherical metric, to a meromorphic function or $\infty$ (see \cite{schiff-1}). Montel's theorem (see \cite{montel-1}) states that if each member of $\mathcal F$ omits three distinct fixed values in $\overline {\mC}$, then $\mathcal F$ is a normal family in $\Omega$. During the last  about hundred years Montel's theorem has undergone various extensions and generalizations. For example, (i) the omitted values are allowed to vary with each member of the family \cite{cara-1}, (ii) the omitted values can be replaced by meromorphic functions \cite{chang-1}  (iii) the omitted values are replaced by mutually avoiding continuous functions \cite{bargmann-1}. 

\medskip

Recall that two nonconstant meromorphic functions $f$ and $g$ defined on $\Omega$ are said to share a set $S$ IM in $\Omega$ provided $\overline E_f(S)=\overline E_g(S)$; where $S$ is a  subset of distinct points in $\overline {\mC}$ and  $$\overline E_f(S):=\bigcup\limits_{a\in S}\{z\in \Omega: f(z)=a\},$$ with each $a$-point in $\overline E_f(S)$ being counted only once.\\

Recently, involving the sharing of values or functions or more generally the sets, various generalizations of Montel's theorem have been obtained (for example, see \cite{chang-1, fang-1, li-1, xu-1, xu-2, xu-3, xu-4}. In particular, M.Fang and W.Hong \cite{fang-1} extended  Montel's Theorem as follows:

\begin{theorem} Let $\mathcal{F}$ be a family of meromorphic functions in a domain $\Omega\subset \mathbb{C}$. If, for each pair of functions $f$ and $g$ in $\mathcal{F}$ share the set $S=\left\{0,1,\infty \right\}$, then the family $\mathcal{F}$ is normal in $\Omega$.
\label{FANGHONG}
\end{theorem}

And J.Chang, M.Fang and L.Zalcman \cite{chang-1} proved the following generalization of Montel's criterion:

\begin{theorem}\label{theorem2} Let $\mathcal{F}$ be a family of meromorphic functions in a domain $\Omega$ and let $a(z)$, $b(z)$ and $c(z)$ be distinct meromorphic functions in $\Omega$, one of which may be $\infty$ identically. If, for all $f\in \mathcal{F}$ and $z\in \Omega$, $f(z)\neq a(z),$ $f(z)\neq b(z)$ and $f(z)\neq c(z)$, then the family $\mathcal{F}$ is normal in $\Omega$.
\end{theorem}

Also, S.Zeng and I.Lahiri \cite{zeng-1} improved the result of Montel by considering shared set of two distinct values and proved the following result:

\begin{theorem}\label{theorem3} Let $\mathcal{F}$ be a family of meromorphic functions in a domain $\Omega$ and $M$ be a positive number and $S=\left\{a, b \right\}$, where $a , b$ are distinct elements of $\overline{\mC}$. Further, suppose that $(i)$ each pair of functions $f,g \in\mathcal{F}$ share the set $S$ in $\Omega$,  $(ii)$ there exists a $c \in \mathbb{C}-\left\{a,b \right\}$ such that for each $f\in \mathcal{F}$, $\left| f'(z)\right |
 \leq M$ whenever $f(z)=c$ in $\Omega$, and $(iii)$ each $f\in \mathcal{F}$ has no simple $b$-points in $\Omega$. Then $\mathcal{F}$ is normal in $\Omega .$
\end{theorem}

In this paper, we extend Theorem \ref{FANGHONG} by replacing the elements of the shared set $S$ by distinct meromorphic functions and hence 
obtain another variation on Montel's theorem.

\begin{theorem}\label{theorem4} Let $\mathcal{F}$ be a family of meromorphic functions in a domain $\Omega$ and let $\psi_1(z)$, $\psi_2(z)$ and $\psi_3(z)$ be distinct meromorphic functions in $\Omega$ such that\\
$(i)$ every $f,g \in\mathcal{F}$ share the set $S:=\left\{\psi_1(z), \psi_2(z), \psi_3(z) \right\}$ in $\Omega$,    \\
$(ii)$ for every $f\in \mathcal{F}$, $f(z_0)\neq \psi_i (z_0)$ whenever $\psi_i(z_0)=\psi_j(z_0)$ for $i,j\in \left\{1,2,3 \right\}(i\neq j)$ and $z_0\in \Omega$.\\
 Then $\mathcal{F}$ is normal in $\Omega$.
\end{theorem}
  
\begin{remark} If for every $f\in \mathcal{F}$, $f(z)\neq \psi_i(z)$ for $i=1,2,3$. Then Theorem \ref{theorem4} reduces to Theorem \ref{theorem2}.
\end{remark}

 As an illustration of Theorem \ref{theorem4}, we have the following example.
	
\begin{example} Consider the family $\mathcal{F}$=$\{$$f_m :m\in\mN\}$, where $$f_m (z)=\frac{e^{z}}{3m}$$ on
the unit disk $\mD$ and let $\psi_1(z)=0$, $\psi_2(z)=e^z$ and $\psi_3(z)= e^z/2$. Clearly, for every $f,g\in\mathcal{F}$, $f$ and $g$ share the set $$S=\left\{0,e^z,\frac{e^z}{2}\right\}$$ of distinct meromorphic functions and  the family $\mathcal{F}$ can easily seen to be normal in $\mD$.
\end{example}

The following example show that the conditions (ii) is essential in Theorem \ref{theorem4}.

\begin{example} Consider the family $\mathcal{F}$=$\{$$f_m :m\in\mN\}$, where $$f_m (z)=2mz$$ on
the unit disk $\mD$ and let $\psi_1(z)=z$, $\psi_2(z)=z/2$ and $\psi_3(z)=z/3$. Clearly, for every $f,g\in\mathcal{F}$, $f$ and $g$ share the set  $$S=\left\{z,\frac{z}{2},\frac{z}{3}\right\}$$ in $\mD$. However, the family $\mathcal{F}$ is not normal in $\mD$. Note that $f_m(0)=\psi_1(0)=\psi_2(0)=\psi_3(0)$, showing that we cannot drop the condition (ii) in Theorem \ref{theorem4}.
\end{example}

Finally, the following example shows that the cardinality of set $S$ in Theorem \ref{theorem4} cannot be reduced.

\begin{example} Consider the family $\mathcal{F}$=$\{$$f_m :m\in\mN\}$, where $$f_m (z)=tan~mz+z$$ on
the unit disk $\mD$ and let $\psi_1(z)=z+i$, $\psi_2(z)=z-i$. Clearly, for every $f,g\in\mathcal{F}$, $f$ and $g$ share the set $$S=\left\{z+i,z-i\right\}$$ and $\psi_1(z)\neq \psi_2(z)$ in $\mD$. However, the family $\mathcal{F}$ is not normal in $\mD$.
\end{example}

Further one can ask what can be said about normality of $\mathcal{F}$ if $f$ is replaced by $f^{(k)}$ in Theorem \ref{theorem4}. In this direction, we prove the following result:

\begin{theorem}\label{theorem5} Let $\mathcal{F}$ be a family of meromorphic functions in a domain $\Omega$, all of whose zeros have multiplicity at least $k+1$, where $k$ is a positive integer. Let $\psi_1(z)$, $\psi_2(z)$ and $\psi_3(z)$ be distinct meromorphic functions in $\Omega$ such that\\
$(i)$ for every $f,g \in\mathcal{F}$, $f^{(k)}$ and $g^{(k)}$ share the set $S:=\left\{\psi_1(z), \psi_2(z), \psi_3(z) \right\}$ in $\Omega$,    \\
$(ii)$ for every $f\in \mathcal{F}$, $f(z_0)\neq \psi_i (z_0)$ whenever $\psi_i(z_0)=\psi_j(z_0)$ for $i,j\in \left\{1,2,3 \right\}(i\neq j)$ and $z_0\in \Omega$.\\
 Then $\mathcal{F}$ is normal in $\Omega$.
\end{theorem}

\begin{example} Consider the family $\mathcal{F}$=$\{$$f_m :m\in\mN\}$, where $$f_m (z)=mz^k$$ on
the unit disk $\mD$ and let $\psi_1(z)=0$, $\psi_2(z)=1/2$ and $\psi_3(z)=1/3$. Clearly, for every $f,g\in\mathcal{F}$, $f^{(k)}$ and $g^{(k)}$ share the set  $$S=\left\{0,\frac{1}{2},\frac{1}{3}\right\}$$ in $\mD$. However, the family $\mathcal{F}$ is not normal in $\mD$. This shows that the condition in Theorem \ref{theorem5} that the zeros of functions in $\mathcal{F}$ have multiplicity at least $k+1$ cannot be dropped.
\end{example}

\section{\textbf{Notations and Lemmas}}
For $z_0\in\mathbb{C}$ and $r>0$, we denote $\mD$ the open unit disk, $D_r(z_0)=\left\{z: \left|z-z_0\right|<r\right\}$ and $D' _r(z_0)=\left\{z: 0<\left|z-z_0\right|<r\right\}$.
 To prove our result, we require the following lemmas.

\begin{lemma}\label{lemma1} (Zalcman's lemma) \cite{schiff-1} Let $\mathcal{F}$ be a family of meromorphic functions in a domain $\Omega$. If $\mathcal{F}$ is not normal at a point $z_0 \in \Omega $, there exist a sequence of points $\left\{z_n\right\}\in \Omega$ with $z_n\rightarrow z_0$, a sequence of positive numbers $\rho_n \rightarrow 0$ and a sequence of functions $f_n \in \mathcal{F}$ such that
 $$g_n(\zeta)=f_n(z_n+\rho_n\zeta)$$ converges locally uniformly with respect to the spherical metric to $g$($\zeta$), where $g(\zeta$) is a non-constant meromorphic function on $\mC$.
\end{lemma}

\begin{lemma}\label{lemma2} \cite{chang-1} Let $\mathcal{F}$ be a family of meromorphic functions in a domain $\mD$ and let a and b be distinct functions holomorphic on $\mD$. Suppose that, for any $f\in\mathcal{F}$ and any $z\in\mD$, $f(z)\neq a(z)$ and $f(z)\neq b(z)$. If $\mathcal{F}$ is normal in $\mD-\left\{0\right\}$, then $\mathcal{F}$ is normal in $\mD$.
\end{lemma}

\begin{lemma}\label{lemma3} \cite{chen-1} Let $\mathcal{F}$ be a family of meromorphic functions in a domain $\Omega$, all of whose zeros have multiplicity at least $k+1$, where $k$ is a positive integer; and let $\mathcal{G}=\left\{f^{(k)}:f\in\mathcal{F}\right\}$. If $\mathcal{G}$ is normal in $\Omega$, then $\mathcal{F}$ is also normal in $\Omega$.
\end{lemma}

\section{\textbf{Proof of Theorems}}

\begin{proof} [\textbf{Proof of Theorem \ref{theorem4}.}]  Since normality is a local property, it is enough to show that $\mathcal{F}$ is normal at each $z_0\in \Omega$. We distinguish the following cases.\\\\
\textbf{Case 1.} $\psi_1 (z_0),\psi_2 (z_0), \psi_3 (z_0)$ are distinct.\\\\
We further consider following subcases.\\\\
\textbf{Case 1.1.} Suppose that there exists $f\in \mathcal{F}$ such that $f(z_0)\neq \psi_i (z_0)$ for $i=1,2,3$. Then we can find a small neighborhood $D_r(z_0)$ in $\Omega$ such that $f(z)\neq \psi_i(z)$ for $i=1,2,3$ in $D_r(z_0)$. By the hypothesis we see that for every $f(z)\in \mathcal{F}$, $f(z)\neq \psi_i(z)$ for $i=1,2,3$ in $D_r(z_0)$. Thus by Theorem \ref{theorem2}, $\mathcal{F}$ is normal at $z_0$.\\\\
\textbf{Case 1.2.} Suppose that there exists $f\in \mathcal{F}$ such that $f(z_0)= \psi_i (z_0)$ for $i=1$ or $2$ or $3$. Without loss of generality we assume $f(z_0)=\psi_2(z_0)=0$ and $\psi_3(z_0)=\infty$.  Since $\psi_i(z_0)\neq \psi_j(z_0)(1\leq i<j\leq 3)$, we can find a small neighborhood $D_r(z_0)$ in $\Omega$ such that $f(z)\neq \psi_i(z)$ for $i=1,2,3$ in $D'_r(z_0)$ and $\psi_1(z)\neq 0,\infty$$(\psi_1\not\equiv 0)$ in $D_r(z_0)$. By the hypothesis we see that for every $f\in \mathcal{F}$, $f(z)\neq \psi_i(z)$ for $i=1,2,3$ in $D'_r(z_0)$. Thus by Theorem \ref{theorem2}, $\mathcal{F}$ is normal in $D'_r(z_0)$. Now we claim that $\mathcal{F}$ is normal at $z_0$.\\
We set  $$ \mathcal{G}:=\left\{g(z)=f(z)-\psi_1(z):f\in\mathcal{F}\right\}.$$
Note that $\mathcal{F}$ is normal if and only if $\mathcal{G}$ is normal. Since $\mathcal{F}$ is normal in $D'_r(z_0)$, so $\mathcal{G}$ is normal in $D'_r(z_0)$. Thus, for a sequence $\left\{g_n\right\}\subset\mathcal{G}$, there exists a subsequence $\left\{ g_{n_k}\right\}$ of $\left\{g_n\right\}$ which converges locally uniformly in $D'_r(z_0)$ to a meromorphic function $h$. We consider the following cases.\\\\
\textbf{Case 1.2.1.} Suppose that $h\equiv \infty$. Then  $$\frac{1}{g_{n_k}(z)}\rightarrow 0   \text{ for }     z\in \partial D_{\frac{r}{2}}(z_0) $$
and since, for any sequence $\left\{f_n\right\}\subset \mathcal{F}$, $f_n(z)$ omits $\psi_1(z)$, we have   $$g_{n_k}(z)\neq 0 \text{ for } z\in D_{\frac{r}{2}}(z_0).$$
Hence there exists $k_0>0$ such that $$\left|\frac{1}{g_{n_k}(z)}\right|\leq M ,~ \forall k\geq k_0,~ z\in \partial D_{\frac{r}{2}}(z_0),$$
where $M>0$ is a constant. Thus by Maximum modulus principle, we conclude that 

$$\left|\frac{1}{g_{n_k}(z)}\right|\leq M ,~ \forall k\geq k_0,~ \forall z\in D_{\frac{r}{2}}(z_0).$$
It follows that $\left\{1 / g_{n_k}(z)\right\}$ converges locally uniformly to 0 in $D_{r/2}(z_0)$ and hence $\left\{ g_{n_k}\right\}$ converges locally uniformly to $h$ in $D_{r/2}(z_0)$. Thus $\mathcal{G}$ is normal at $z_0$.\\\\
\textbf{Case 1.2.2.} Suppose that $h\not\equiv \infty$. Then again there exists an index $k_0 > 0$ such that $$\left|g_{n_k}(z)\right|\leq M,~\forall k\geq k_0,~ z\in \partial D_{\frac{r}{2}}(z_0) $$ and since, for any sequence $\left\{f_n\right\}\subset \mathcal{F}$, $f_n(z)\neq \infty$ and $\psi_1(z)\neq \infty$, we have $$g_{n_k}(z)\neq \infty \text{~for~} z\in D_{\frac{r}{2}}(z_0),$$
where $M>0$ is a constant. Thus by Maximum modulus principle, we conclude that $$\left|g_{n_k}(z)\right|\leq M,~\forall k\geq k_0,~ \forall z\in D_{\frac{r}{2}}(z_0).$$
It follows that $\left\{g_{n_k}(z)\right\}$ converges locally uniformly to $h$ in $D_{r/2}(z_0)$. Hence there exists a subsequence of $\left\{g_{n}(z)\right\}$ which converges locally uniformly to $h$ in $D_{r/2}(z_0)$. Therefore $\mathcal{G}$ is normal at $z_0$.

Thus, $\mathcal {F}$ is normal at $z_0$.\\\\
\textbf{Case 2.} Exactly two of $\psi_1(z_0), \psi_2(z_0), \psi_3(z_0)$ are equal.\\\\
Without loss generality we assume that $\psi_1(z_0)=\psi_2(z_0)$ and $\psi_3(z_0) \neq \psi_1(z_0),\psi_2(z_0)$. Then, by hypothesis, for every $f\in \mathcal{F}$, $f(z_0)\neq \psi_i(z_0)$ for $i=1,2$. We consider the following two subcases.\\\\
\textbf{Case 2.1.} $\psi_1(z_0)$ is finite.\\\\
Then we can find a small neighborhood $D_r(z_0)$ in $\Omega$ such that $\psi_i(z)\neq \psi_j(z)$ $(1\leq i<j\leq 3)$ in $D'_r(z_0)$. Thus by Case 1, $\mathcal{F}$ is normal in $D'_r(z_0)$. Now we turn to show that $\mathcal{F}$ is also normal at $z_0$. Since $\psi_1(z_0)$ and $\psi_2(z_0)$ are finite and for every $f\in \mathcal{F}$, $f(z_0)\neq \psi_i(z_0)$ for $i=1,2$, we can find that for every $f\in \mathcal{F}$, $f(z)\neq \psi_i(z)$ for $i=1,2$ and $\psi_1(z),\psi_2(z)$ are holomorphic in $D_r(z_0)$. Then by Lemma \ref{lemma2}, $\mathcal{F}$ is normal at $z_0$.\\\\
\textbf{Case 2.2.} $\psi_1(z_0)$ is infinite.\\\\
Then we can find a small neighborhood $D_r(z_0)$ in $\Omega$ such that $\psi_1(z)$ and $\psi_2(z)$ are holomorphic in $D'_r(z_0)$ and $\psi_i(z)\neq 0$ for $i=1,2$ in $D_r(z_0)$.\\
 We set
$$\mathcal{G}=\left\{g=\frac{1}{f} :f\in \mathcal{F}\right\}$$
and

 $$\phi_i(z)=\frac{1}{\psi_i(z)}~~\text{for}~~ i=1,2,3.$$
Then $\phi_1(z_0)=\phi_2(z_0)=0$ and $\phi_3(z_0)\neq 0$. Thus, as in Case 2.1, we can prove that $\mathcal{G}$ is normal at $z_0$ and hence $\mathcal{F}$ is normal at $z_0$.\\\\
\textbf{Case 3.} $\psi_1(z_0)=\psi_2(z_0)=\psi_3(z_0)$.\\\\
By the hypothesis, we have for every $f\in\mathcal{F}$, $f(z_0)\neq \psi_i(z_0)$ for $i=1,2,3$. Then we can find a small neighborhood $D_r(z_0)$ in $\Omega$ such that for every $f\in\mathcal{F}$, $f(z)\neq \psi_i(z)$ for $i=1,2,3$ in $D_r(z_0)$. Hence by Theorem \ref{theorem2}, $\mathcal{F}$ is normal at $z_0$.

This completes the proof of Theorem \ref{theorem4}.
\end{proof}

\begin{proof} [\textbf{Proof of Theorem \ref{theorem5}.}] Since normality is a local property, it is enough to show that $\mathcal{F}$ is normal at each $z_0\in \Omega$. We distinguish the following cases.\\\\
\textbf{Case 1.} $\psi_i (z_0)\neq \psi_j (z_0)$ for $1\leq i < j \leq 3$.\\\\
Then we can find a small neighborhood $D_r(z_0)$ in $\Omega$ such that $\psi_i(z)\neq \psi_j(z)$ $(1\leq i<j\leq 3)$ in $D_r(z_0)$. Since, for every $f,g \in\mathcal{F}$, $f^{(k)}$ and $g^{(k)}$ share the set $S=\left\{\psi_1(z), \psi_2(z), \psi_3(z) \right\}$. Thus by Theorem \ref{theorem4}, $\left\{f^{(k)}:f\in \mathcal{F}\right\}$ is normal in $D_r(z_0)$, so by Lemma \ref{lemma3}, $\mathcal{F}$ is normal in $D_r(z_0)$. Thus $\mathcal{F}$ is normal at $z_0$.\\\\
\textbf{Case 2.} $\psi_i (z_0)= \psi_j (z_0)$ for $1\leq i < j \leq 3$.\\\\
Proceeding in the same way as in Case 2 and Case 3 of Theorem \ref{theorem4}, we conclude that $\mathcal{F}$ is normal at $z_0$. 

This completes the proof of Theorem \ref{theorem5}.
\end{proof}

\bibliographystyle{plain}

\end {document}